\documentclass[12pt]{amsart}
\usepackage{geometry}
\usepackage{cite}
\usepackage{hyperref}
\usepackage{amssymb}
\usepackage{mathrsfs}

\newtheorem{theorem}{Theoreme}[section]

\newtheorem{question}{Question}

\newtheorem{definition}[theorem]{Definition}

\newcommand{\dsum}{\mathop{\oplus }}

\begin{document}
\vspace*{2cm}

\title{Vector Fourier analysis on compact groups and Assiamoua spaces}

\author{Yaogan Mensah}
\address{Department of Mathematics, University of Lom\'e, POB 1515 Lom\'e 1,Togo}
\address{and  International Chaire in Mathematical Physics and Applications, University  of Abomey-Calavi, Benin}
\email{mensahyaogan2@gmail.com, ymensah@univ-lome.tg}

\dedicatory{In memory of Professor VSK Assiamoua }

\begin{abstract}
This paper shows how a family of function spaces, coined as  Assiamoua spaces, plays a fundamental role in the Fourier analysis of vector-valued functions on  compact groups. These spaces make it possible to transcribe the classical results of Fourier analysis in the framework of analysis of vector-valued functions and vector measures. The construction of Sobolev spaces of vector-valued functions on compact groups rests heavily on the members of the aforementioned family.  
\end{abstract}
\maketitle
\section{Introduction}
The Fourier analysis has undergone many generalizations in various directions since the original idea of J. Fourier \cite{Fourier}. The extension of this analysis to complex-valued functions on a compact group has become a well documented classical theory. While Fourier used trigonometric series, today known as Fourier series, in the modeling of the heat phenomenon, his theory has been extended to abelian groups, compact groups, locally compact groups, then to other more complex topological and algebraic structures such as locally compact hypergroups. Each extension is accompanied by the delicate construction of a family of function spaces in which the Fourier transforms of the initial functions reside. For instance, in the case of complex-valued functions on the one-dimensional torus or on $\mathbb{R}^n$, Lebesgue spaces and their discrete analogues  made it possible to report the results in terms of Riemann-Lebesgue lemma, Fourier inversion formula, Plancherel theorem, etc. For more details, we refer the books  \cite{Folland, Grafakos}. 

In the noncommutative setting, for instance for compact non-necessarily abelian groups or  for more general  locally compact groups, the Fourier analysis of complex functions appeals  classes of bounded operators  on the representation spaces. The case of compact groups is described in the present paper. For more details on noncommutative Fourier analysis, we refer the books  \cite{Folland2, Hewitt}.  

The situation becomes more complicated when we consider vector-valued functions. The study of the Fourier analysis of vector-valued functions and vector measures was initiated by Assiamoua in \cite{Assiamoua1} and developed by Assiamoua and Olubummo in\cite{Assiamoua2}. A family of functions spaces, coined as Assiamoua spaces in the present paper, was introduced in order to conduct effectively the study.  These spaces are quite useful and they recently made it possible to obtain an extension of the Fourier-Stieltjes transform from compact groups to locally compact groups through  an induced representation \cite{Akakpo}.

The aim of this paper is to highlight these function spaces by emphasizing their roles  and the potential of what they can still bring to the development of abstract harmonic analysis.

The rest of the paper is organized  as follows. In  Section \ref{Basics of  representation theory},we recall basic definitions and facts in group representation theory. Section \ref{Fourier analysis on compact groups} provides a summary of the Fourier analysis of compact groups. Section \ref{Assiamoua spaces} is devoted to the  Fourier analysis of vector-valued functions  on compact groups and  the highlighting of Assiamoua spaces. In Section \ref{Tensorization}, a tensor product approach is used to describe the Assiamoua spaces in a way equivalent  to the previous construction. In Section \ref{Topological  properties}, some topological properties of the Assiamoua spaces are provided. In Section \ref{Sobolev}, we describe a Sobolev space whose construction  is based on an Assiamoua space. Finally, Section \ref{Ideas} collects some ideas and questions  for further research.   

\section{Basics of  representation theory}\label{Basics of  representation theory}
In this section, we recall some basic notions in group representation theory in order to make the paper more understandable.   

Let $G$ be a topological compact. 
 A unitary representation $\sigma$ of $G$  on a complex  Hilbert space $H_\sigma$ is a  continuous group homomorphism  $\sigma : G\longrightarrow \mathcal{U}(H_\sigma)$ where  $\mathcal{U}(H_\sigma)$ denotes the group of unitary operators on $H_\sigma$. 
The continuity is to be understood as the mapping $$G\longrightarrow H_\sigma, x\longmapsto \sigma (x)\xi$$ is continuous for every $\xi\in H_\sigma$. 
  The Hilbert space $H_\sigma$ is called the representation space of $\sigma$ and the dimension of $H_\sigma$  is called the dimension of the representation $\sigma$ and it is usually denoted by $d_\sigma$. 
  
  A subspace $M$ of $H_\sigma$ is said to be invariant by $\sigma$ if $\forall x\in G, \forall \xi \in M, \sigma (x)\xi\in M$. 
 A  representation $\sigma$ of $G$ on  $H_\sigma$ is said to be  irreducible if there is no proper closed subspace $M$ of $H_\sigma$ which is  invariant by $\sigma$.  
 
It is well known that  the dimension of any  unitary  irreducible  representation of a compact group  is  of finite dimension \cite{Hewitt}. 

Two unitary representions $\sigma_i,\, i=1,2$ of $G$ on $H_{\sigma_i},\,i=1,2$ are said to be unitary equivalent if there exists a unitary  operator $S: H_{\sigma_1}\longrightarrow H_{\sigma_2}$
such that $\forall x\in G, \,\sigma_2(x)S=S\sigma_1(x)$.

Throughout the paper, $\widehat{G}$ stands for  the  set of equivalent classes  of unitary irreducible  representations of $G$. It is called the unitary  dual of $G$. If  $G$ is compact, then  $\widehat{G}$  is discrete. 

\section{Fourier analysis on compact groups}\label{Fourier analysis on compact groups}
In this section, fundamental concepts and results of the Fourier analysis of complex-valued functions  on compact groups are recaped. More details can be find in many abstract harmonic analysis books. However, we refer the preliminary notes of the paper \cite[Section 2]{Kumar}. Results in this section are generalized in Section \ref{Assiamoua spaces}. 

Let $G$ be a compact group (assumed Hausdorff once for all) equipped with its normalized Haar measure. Let us denote by $L^1(G)$ the space of complex-valued integrable functions on $G$ and by $L^2(G)$ the space of complex-valued square integrable functions on $G$. The Fourier transform of $f\in  L^1(G)$ is defined by 

$$\widehat{f}(\sigma)=\int_Gf(x)\sigma (x)^*dx,\, \sigma \in \widehat{G},$$
where  $\sigma(x)^*$ stands for the Hilbert adjoint of the operator $\sigma(x)$. The object $\widehat{f}(\sigma)$ is a bounded linear operator on $H_\sigma$.
Then, $\widehat{f}=\left(\widehat{f}(\sigma)\right)_{\sigma \in \widehat{G}}$ is a family of  bounded linear operators. Let us denote  by  $B(H_\sigma)$  the space of bounded linear operators on $H_\sigma$ endowed with the  operator norm. To express results about the Fourier transform of complex-valued functions on compact groups, the following spaces are considered.

\begin{enumerate}
\item $\ell^p-\dsum\limits_{\sigma\in \widehat{G}}B(H_\sigma)=\left\lbrace (T_\sigma)\in \prod\limits_{\sigma\in \widehat{G}}B(H_\sigma) : \sum\limits_{\sigma\in \widehat{G}}d_\sigma \|T_\sigma\|_{B(H_\sigma)}^p<\infty\right\rbrace,\,1\leqslant p<\infty$, endowed with the norm 
$$ \|(T_\sigma)\|_{\ell^p-\dsum\limits_{\sigma\in \widehat{G}}B(H_\sigma)}=\left(\sum\limits_{\sigma\in \widehat{G}}d_\sigma \|T_\sigma\|_{B(H_\sigma)}^p\right)^{\frac{1}{p}},$$ 

\item $\ell^\infty-\dsum\limits_{\sigma\in \widehat{G}}B(H_\sigma)=\left\lbrace (T_\sigma)\in \prod\limits_{\sigma\in \widehat{G}}B(H_\sigma) : \sup\limits_{\sigma\in \widehat{G}}\|T_\sigma\|_{B(H_\sigma)}<\infty\right\rbrace$, endowed with the norm 
$$\|(T_\sigma)\|_{\ell^\infty-\dsum\limits_{\sigma\in \widehat{G}}B(H_\sigma)}=\sup\limits_{\sigma\in \widehat{G}}\|T_\sigma\|_{B(H_\sigma)},$$

\item  $c_0-\dsum\limits_{\sigma\in \widehat{G}}B(H_\sigma)=\left\lbrace (T_\sigma)\in \ell^\infty-\dsum\limits_{\sigma\in \widehat{G}}B(H_\sigma): T_\sigma \longrightarrow 0 \mbox{ as } \sigma\longrightarrow \infty\right\rbrace$. 
 \end{enumerate}

The spaces thus constructed make it possible to efficiently transcribe the results of harmonic analysis on compact groups. Let us mention some well known results. 

\begin{theorem}
Let $G$ be a compact group. If $f\in L^1(G)$, then $\widehat{f} \in \ell^\infty-\dsum\limits_{\sigma\in \widehat{G}}B(H_\sigma)$. Moreover, 
$$\|\widehat{f}\|_{\ell^\infty-\dsum\limits_{\sigma\in \widehat{G}}B(H_\sigma)}\leqslant \|f\|_{L^1(G)}.$$
\end{theorem}
\begin{theorem}Let $G$ be a compact group. If $f\in L^1(G)$, then $\widehat{f} \in c_0-\dsum\limits_{\sigma\in \widehat{G}}B(H_\sigma)$. 
\end{theorem}
The above theorem is known as the Riemann-Lebesgue lemma.

\begin{theorem}
Let $G$ be a compact group.
 The mapping $f\longrightarrow \widehat{f}$ is an isomorphism  from $L^2(G)$ onto $\ell^2-\dsum\limits_{\sigma\in \widehat{G}}B(H_\sigma)$. Moreover, if $f\in L^2(G)$, then
 $$\|\widehat{f}\|_{\ell^2-\dsum\limits_{\sigma\in \widehat{G}}B(H_\sigma)}=\|f\|_{L^2(G)}.$$
\end{theorem}
The above theorem is known as the Plancherel theorem.

\section{Fourier analysis of vector-valued functions on compact groups and Assiamoua spaces}\label{Assiamoua spaces}

Let $G$ be a compact group and let $A$ be a complex  Banach space. Let $L^1(G,A)$ denotes the set of $A$-valued functions on $G$ which are integrable in the Bochner sense. The set $L^1(G,A)$ is endowed with the norm

$$\|f\|_{L^1(G,A)}=\int_G\|f(x)\|_A dx.$$

In 1989, in order to study the multipliers of $L^1(G, A)$ \cite{Assiamoua1}, Assiamoua  introduced the Fourier-Stieltjes transform of vector measures on compact groups and consequently the Fourier transform of  vector-valued functions on compact groups. 

Let $f\in L^1(G, A)$. The Fourier transform of  $f$ has been interpreted for the first time  by Assiamoua as the family $\widehat{f}=(\widehat{f}(\sigma))_{\sigma\in\widehat{G}}$ of sesquilinear maps $\widehat{f}(\sigma): H_\sigma\times H_\sigma \longrightarrow A $ defined by
\begin{equation}\label{Assiamoua's formula}
\widehat{f}(\sigma)(\xi,\eta)=\int_G\langle\sigma(x)^*\xi,\eta\rangle_{H_\sigma}f(x)dx, \, \xi,\eta\in H_\sigma.
\end{equation}
The formula (\ref{Assiamoua's formula}) is expressed here using our own notations (see \cite[Section 4]{Assiamoua1}).  From now on, the formula (\ref{Assiamoua's formula}) will be called  Assiamoua's formula. The systematic study of the main properties of this transformation was made by Assiamoua and Olubummo in \cite{Assiamoua2}. To efficiently transcribe the results  for the aformentioned study, some function spaces were introduced : the analogue of the spaces of Section \ref{Fourier analysis on compact groups}. Let us briefly describe their construction. 

For each $\sigma\in\widehat{G}$, fix once for all an orthonormal basis $(\xi_1^\sigma, \cdots ,\xi_{d_\sigma}^\sigma )$ of $H_\sigma$. 
Denote by  
$\mathscr{S}(H_\sigma\times H_\sigma, A)$  the set of sesquilinear mappings from $H_\sigma\times H_\sigma$ into $A$. The set $\mathscr{S}(H_\sigma\times H_\sigma, A)$ is endowed with the norm 
$$\|\psi\|_{\mathscr{S}(H_\sigma\times H_\sigma, A)}=\sup\{\|\psi (\xi,\eta)\|_A: \|\xi\|\leqslant 1, \|\eta\|\leqslant 1\}.$$
Set 
$$\mathscr{S}(\widehat{G},A)=\prod\limits_{\sigma\in\widehat{G}}\mathscr{S}(H_\sigma\times H_\sigma, A).$$
Consider the following spaces :
\begin{enumerate}
\item $\mathscr{S}_p(\widehat{G},A)=\left\lbrace \varphi \in \mathscr{S}(\widehat{G},A) : \sum\limits_{\sigma \in \widehat{G}} d_\sigma \sum\limits_{i=1}^{d_\sigma}\sum\limits_{j=1}^{d_\sigma}\|\varphi (\sigma)(\xi_j^\sigma,\xi_i^\sigma )\|_A^p<\infty \right\rbrace,\,  1\leqslant p<\infty$ equipped with the norm 
$$\|\varphi\|_{\mathscr{S}_p(\widehat{G},A)}=\left(  \sum\limits_{\sigma \in \widehat{G}} d_\sigma \sum\limits_{i=1}^{d_\sigma}\sum\limits_{j=1}^{d_\sigma}\|\varphi (\sigma)(\xi_j^\sigma,\xi_i^\sigma )\|_A^p \right)^{\frac{1}{p}},$$

\item $\mathscr{S}_\infty(\widehat{G},A)=\left\lbrace \varphi \in \mathscr{S}(\widehat{G},A) : \sup\limits_{\sigma\in\widehat{G}}\|\varphi (\sigma)\|_{\mathscr{S}(H_\sigma\times H_\sigma, A)}<\infty\right\rbrace$ equipped with the norm
$$ \|\varphi\|_{\mathscr{S}_\infty(\widehat{G},A)}=  \sup\limits_{\sigma\in\widehat{G}}\|\varphi (\sigma)\|_{\mathscr{S}(H_\sigma\times H_\sigma, A)},$$
\item  $\mathscr{S}_0(\widehat{G},A)=\left\lbrace  \varphi \in \mathscr{S}_\infty(\widehat{G},A) : \forall \varepsilon >0,\{\sigma\in \widehat{G} : \|\varphi (\sigma)\|_{\mathscr{S}(H_\sigma\times H_\sigma, A)}>\varepsilon\} \mbox{ is finite } \right\rbrace.$
\end{enumerate}

\begin{definition}\label{Definition of Assiamoua spaces}
The spaces  $\mathscr{S}_p(\widehat{G},A), \, 1\leqslant p\leqslant \infty$ are called Assiamoua spaces. 
\end{definition}

Using the Assiamoua spaces, some properties of the  Fourier transform read as follows \cite{Assiamoua2}.

\begin{theorem} Let $G$ be a compact group and $A$ a complex Banach space.
The mapping $f\mapsto \widehat{f}$ from $L^1(G,A)$ into $\mathscr{S}_\infty(\widehat{G},A)$ is linear, injective and continuous. 
\end{theorem}
The following theorem contains the  analogue of the  Riemann-Lebesgue lemma. 
\begin{theorem}Let $G$ be a compact group and $A$ a complex Banach space.
If $f\in L^1(G, A)$, then $\widehat{f}\in \mathscr{S}_0(\widehat{G},A)$. Moreover, $\widehat{L^1(G, A)}$ is dense in $\mathscr{S}_0(\widehat{G},A)$ where $\widehat{L^1(G, A)}$ is the image of $L^1(G, A)$ by the Fourier transform.
\end{theorem}

\begin{theorem}\label{Plancherel}
Let $G$ be a compact group. The Fourier transformation $f\mapsto \widehat{f}$ is an isometry from $L^2(G,A)$ onto $\mathscr{S}_2(\widehat{G},A)$. Moreover, if $f\in L^2(G,A)$ then 
$$ f(x)= \sum\limits_{\sigma \in \widehat{G}} d_\sigma \sum\limits_{i=1}^{d_\sigma}\sum\limits_{j=1}^{d_\sigma}\widehat{f}(\sigma)\langle \sigma (x)\xi_j^\sigma,\xi_i^\sigma\rangle,\, x\in G.$$
\end{theorem}

\section{Tensor product approach to Assiamoua spaces }\label{Tensorization}
Using the  tensor product theory, it is possible to construct  the Assiamoua spaces of  Section \ref{Assiamoua spaces}  in an equivalent way. Denote by $\overline{H}_\sigma$ the conjugate of the Hilbert space $H_\sigma$. The tensor product space $H_\sigma \otimes \overline{H}_\sigma$ is equipped with the projective tensor norm. Since it is  of finite dimension (its dimension is $d_\sigma^2$) then it  is complete. 

One may interpret the Fourier  transform as the family $(\widehat{f}(\sigma))_{\sigma \in \widehat{G}}$ of linear operators from $H_\sigma \otimes \overline{H}_\sigma$ into $A$ defined by (see \cite{Mensah}) :
\begin{equation}
\widehat{f}(\sigma)(\xi\otimes \eta)=\int_G\langle\sigma(x)^*\xi,\eta\rangle_{H_\sigma}f(x)dx, \, \xi\in H_\sigma,\eta\in\overline{H}_\sigma.
\end{equation} 
The Assiamoua spaces as defined in Section \ref{Assiamoua spaces} are slightly modified. 

Denote by $\mathfrak{B}(H_\sigma \otimes \overline{H}_\sigma, A)$ the space of bounded linear operators from  $H_\sigma \otimes \overline{H}_\sigma$ into $A$ equipped with the norm

$$\|\psi\|_{\mathfrak{B}(H_\sigma \otimes \overline{H}_\sigma, A)}=\sup\{\|\psi(\xi\otimes \eta)\|_A : \|\xi\|_{H_\sigma}\leqslant 1, \|\eta\|_{\overline{H}_\sigma}\leqslant 1\}.$$

Now, set 
$$\mathfrak{B}(\widehat{G},A)=\prod\limits_{\sigma \in\widehat{G}}\mathfrak{B}(H_\sigma \otimes \overline{H}_\sigma, A).$$

The space $\mathscr{S}_p(\widehat{G},A)$ becomes $$\mathfrak{B}_p(\widehat{G},A)=\left\lbrace \varphi \in \mathfrak{B}(\widehat{G},A) : \sum\limits_{\sigma \in \widehat{G}} d_\sigma \sum\limits_{i=1}^{d_\sigma}\sum\limits_{j=1}^{d_\sigma}\|\varphi (\sigma)(\xi_j^\sigma\otimes \xi_i^\sigma )\|_A^p<\infty \right\rbrace,\,  1\leqslant p<\infty$$ equipped with the norm 

$$ \|\varphi\|_{\mathfrak{B}_p(\widehat{G},A)}=\left(\sum\limits_{\sigma \in \widehat{G}} d_\sigma \sum\limits_{i=1}^{d_\sigma}\sum\limits_{j=1}^{d_\sigma}\|\varphi (\sigma)(\xi_j^\sigma\otimes \xi_i^\sigma )\|_A^p\right)^{\frac{1}{p}}$$

and $\mathscr{S}_\infty(\widehat{G},A)$ becomes

$$\mathfrak{B}_\infty(\widehat{G},A)=\left\lbrace \varphi \in \mathfrak{B}(\widehat{G},A) : \sup\limits_{\sigma\in\widehat{G}}\|\varphi (\sigma)\|_{\mathfrak{B}(H_\sigma\otimes \overline{H}_\sigma, A)}<\infty\right\rbrace$$ equipped with the norm
$$ \|\varphi\|_{\mathfrak{B}_\infty(\widehat{G},A)}=  \sup\limits_{\sigma\in\widehat{G}}\|\varphi (\sigma)\|_{\mathfrak{B}(H_\sigma\otimes \overline{H}_\sigma, A)}.$$

One can see that the correspondance $\mathscr{S}_p(\widehat{G},A) \leftrightarrow \mathfrak{B}_p(\widehat{G},A)$ is one-to-one. In fact, if we denote by $B(H_\sigma\times \overline{H}_\sigma, A)$ the space of bilinear maps from $H_\sigma\times \overline{H}_\sigma$ into $A$, then we have the following identification  by the effect of linearization of bilinear maps :
$$ \mathscr{S}(H_\sigma \times H_\sigma, A)\simeq  B(H_\sigma\times \overline{H}_\sigma, A)\simeq \mathfrak{B}(H_\sigma\otimes \overline{H}_\sigma, A).$$
The following definition is equivalent to Definition \ref{Definition of Assiamoua spaces}.  
\begin{definition}
We call Assiamoua spaces the spaces $\mathfrak{B}_p(\widehat{G},A),\, 1\leqslant p \leqslant \infty$. 
\end{definition}
The well known properties of the Fourier analysis of vector-valued functions on compact groups  can be transcribed using  the spaces $\mathfrak{B}_p(\widehat{G},A)$. For instance, The Fourier transformation $f\mapsto\widehat{f}$ is an isometry from $L^2(G,A)$ onto $\mathfrak{B}_2(\widehat{G},A)$.

\section{Topological properties of Assiamoua spaces}\label{Topological properties}
The Assiamoua spaces $\mathscr{S}_p(\widehat{G},A)$ were studied from the topological point of view. It had been found out that they have interesting topological  properties. The detailed proofs can be found in \cite{Assiamoua2, MensahAssiamoua1}. Let us mention some of the results. 
\begin{theorem}
The spaces $\mathscr{S}_p(\widehat{G},A),\, 1\leqslant p \leqslant \infty$,  are Banach spaces.
\end{theorem}
 The case $p=\infty$ was studied in \cite{Assiamoua2} and the case $1\leqslant p < \infty$ was studied in \cite{MensahAssiamoua1}.  
 
 \begin{theorem}\cite{Assiamoua2}
 If $A$ is a  Hilbert space, then so is  $\mathscr{S}_2(\widehat{G},A)$. 
 \end{theorem}
 In this case, the inner product in $\mathscr{S}_2(\widehat{G},A)$ is defined as follows : for $\varphi_1, \varphi_2\in \mathscr{S}_2(\widehat{G},A)$, 
\begin{equation}
\langle \varphi_1, \varphi_2\rangle_{\mathscr{S}_2(\widehat{G},A)}=\sum\limits_{\sigma\in \widehat{G}}d_\sigma\sum\limits_{i=1}^{d\sigma}\sum\limits_{j=1}^{d\sigma}\langle \varphi_1(\sigma)(\xi_j^\sigma,\xi_i^\sigma ), \varphi_2(\sigma)(\xi_j^\sigma,\xi_i^\sigma )\rangle_A. 
\end{equation} 
 
\begin{theorem}\label{dual}
The topological dual  of the space $\mathscr{S}_p(\widehat{G},A),\,1\leqslant p < \infty$ is $\mathscr{S}_q(\widehat{G},A^*)$ where $q$ is such that $\displaystyle\frac{1}{p}+\frac{1}{q}=1$ and $A^*$ is the topological dual  of $A$.
\end{theorem}
 The case $p=1$ was proved in \cite{MensahAssiamoua2} when studying the  Fourier algebra of vector-valued functions on compact groups.  The case $1<p<\infty$ was proved in \cite{MensahAssiamoua1}. In both cases, the results are achieved  by showing  that the following mapping is a linear surjective isometry : 
 
 $$T: \mathscr{S}_q(\widehat{G},A^*) \longrightarrow \left(\mathscr{S}_p(\widehat{G},A)\right)^*$$ 
defined by 

$$\langle T\theta,\varphi \rangle=\sum\limits_{\sigma\in\widehat{G}}d_\sigma\sum\limits_{i=1}^{d_\sigma}\sum\limits_{j=1}^{d_\sigma}\langle \theta (\sigma)(\xi_j^\sigma,\xi_i^\sigma),\varphi (\sigma)(\xi_j^\sigma,\xi_i^\sigma) \rangle,$$
$\theta \in \mathscr{S}_q(\widehat{G},A^*)$, $\varphi\in \mathscr{S}_p(\widehat{G},A)$ and the brackets $\langle \cdot, \cdot\rangle$ of the left and right sides are duality brackets.

\section{The Sobolev spaces that Assiamoua could have constructed}\label{Sobolev}
When we were writing the manuscript \cite{Mensah2} on Sobolev spaces of Bessel potential type consisting of vector-valued functions on a compact group, we told ourself that Assiamoua could have written on the subject in 1989 because he had all the necessary ingredients at his disposal in the paper \cite{Assiamoua2}  : his elegant  interpretation of the Fourier transform for vector-valued functions, a Fourier inversion formula and an analogue of the Plancherel theorem. All he had to do was to put them together to  complete the puzzle.  Unfortunately, he did not do so probably because he was not interested in Sobolev spaces at that time. Moreover, as far as we know, it was in 2014 that was published the first article dealing with Sobolev spaces defined on abstract locally compact abelian groups using the Fourier transformation  without taking into account the weak derivation\cite{Gorka}. 

The Sobolev spaces $H_\gamma^s(G,A)$ constructed in  \cite{Mensah2} rests on the Assiamoua space $\mathscr{S}_2(\widehat{G},A)$. Let $s\in [0,\infty)$ and consider the map $\gamma : \widehat{G}\longrightarrow [0,\infty)$. Then, the Sobolev space $H_\gamma^s(G,A)$ is defined as 
$$H_\gamma^s(G,A)=\left\lbrace f\in L^2(G,A) : (1+\gamma^2)^{\frac{s}{2}}\widehat{f}\in  \mathscr{S}_2(\widehat{G},A)\right\rbrace,$$
endowed with the norm 

$$\|f\|_{H_\gamma^s(G,A)}=\left(\sum\limits_{\sigma\in \widehat{G}}d_\sigma(1+\gamma (\sigma)^2)^s\sum\limits_{i=1}^{d_\sigma}\sum\limits_{j=1}^{d_\sigma}\|\widehat{f}(\sigma)(\xi_j^\sigma,\xi_i^\sigma)\|_A^2\right)^{\frac{1}{2}}.
$$
Using properties of the Fourier transform (mainly the inversion formula and the Plancherel type theorem) together with those of the space $\mathscr{S}_2(\widehat{G},A)$, we obtain results,  some of which are recorded in the following theorem \cite{Mensah2}:
\begin{theorem}
Let $G$ be a compact group. Let $A$ be a  Banach space.
\begin{enumerate}
\item The  space  $H_\gamma^s(G,A)$ is a  Banach space. 
\item The Sobolev space $H_\gamma^s(G,A)$ is continuously embedded  in the space $L^2(G,A)$ with $$\|f\|_{L^2(G,A)}\leqslant \|f\|_{H_\gamma^s(G,A)}.$$
\item  If $\sum\limits_{\sigma\in\widehat{G}}\displaystyle\frac{d_\sigma ^3}{(1+\gamma(\sigma)^2)^s}<\infty$, then the Sobolev space $H_\gamma^s(G,A)$ is continuously embedded in the space $\mathcal{C}(G,A)$, where  $\mathcal{C}(G,A)$ is the space of $A$-valued continuous functions on $G$.
\end{enumerate}
\end{theorem}

\section{Some ideas to move further}\label{Ideas}
In this section, we discuss  some ideas in a questioning form  for future research related to  the Assiamoua spaces $\mathscr{S}_p(\widehat{G},A), 1\leqslant p\leqslant \infty$. 
\begin{enumerate}
\item The space $\mathscr{S}_p(\widehat{G},A)$  is a Banach space. However, it can possess other interesting properties which can be investigated. 
\begin{question}
Does it have geometric properties such as the Radon-Nikodym property, the Dunford-Pettis property or the Schur property?
\end{question}

\item The Pego theorem characterizes the precompact subsets of $L^2(\mathbb{R}^n)$ through the Fourier transform. The analogue
of the Pego theorem on compact groups was proved in \cite{Kumar} via the Fourier analysis of complex-valued functions on compact groups.
\begin{question}
Can we prove the analogue of the Pego theorem for $L^2(G,A)$ using the Fourier transform of vector-valued functions and some properties of the spaces $\mathscr{S}_p(\widehat{G},A)$?
\end{question}
\item The Rellich-Kondrachov theorem provides conditions under which a Sobolev space is  compactly embedded  in a Lebesgue space. 
\begin{question}
Can we investigate  the analogue of the  Rellich-Kondrachov theorem for the Sobolev space $H_\gamma^s(G,A)$? 
\end{question}
\item We may define and study the Sobolev spaces of Section \ref{Sobolev} in a more general setting. For instance, we can  consider the Sobolev spaces $H_\gamma^{s,p}(G,A),\, 1\leqslant p\leqslant \infty$ defined by 
$$H_\gamma^{s,p}(G,A)=\left\lbrace f\in L^2(G,A) : (1+\gamma^2)^{\frac{s}{2}}\widehat{f}\in  \mathscr{S}_p(\widehat{G},A)\right\rbrace, \, 1\leqslant p\leqslant \infty.$$
 \begin{question}
 What are the main properties of the Sobolev spaces $H_\gamma^{s,p}(G,A)$ and what are their connections with other function spaces? 
\end{question}
\end{enumerate}

\section{Conclusion}
In this article, we have briefly described the Fourier analysis of complex-valued functions and vector-valued functions on compacts groups. We have highlighted a family of function spaces coined as Assiamoua spaces which play a fundamental role in the analysis of vector-valued functions and we have exhibited some of their major properties.  Finally, we have provided a list of questions which may interest some reseachers in the field of harmonic analysis or functional analysis in general.

\end{document}